\documentclass[11pt,letter]{article}
\usepackage{amsmath}
\usepackage{amsfonts}
\usepackage{amssymb}
\usepackage{graphicx}
\usepackage{enumerate}
\usepackage[ruled]{algorithm2e}
\usepackage{tabularx,booktabs}
\usepackage{float}

\newtheorem{theorem}{Theorem}
\newtheorem{proposition}[theorem]{Proposition}

\newcommand{\QED}{\hfill$\square$}

\title {
    \bf {A note on 'A New Approach To Compute Wiener Index'}
}

\author
{
{\large \sc Aleksandar Ili\' c } \\
{\em \normalsize Facebook Inc, Menlo Park, CA, USA} \\
{\normalsize e-mail: { \tt aleksandari@gmail.com }}
}

\begin{document}

\maketitle

\begin{abstract}
In this note, we discuss the method explained in the recent paper [P. Manuel, I. Rajasingh, B. Rajan, R. Sundara Rajan, \emph{A New Approach To Compute Wiener Index}, Journal of Computational and Theoretical Nanoscience {\bf 10}, (2013) 1515--1521.] for computing the Wiener index of special chemical graphs. The method is actually already well-known and equivalent to the 'cut method' introduced in 1995 by Klav\v zar, Gutman and Mohar, and used in multiple papers for computing various distance based graph invariants.
\end{abstract}

\vspace{0.2cm}

\section{Preliminaries}

Let $G = (V, E)$ be a simple connected graph. The distance $d_G(u, v)$ between two vertices $u$ and $v$
is defined as the number of edges on a  path connecting $u$ and $v$.

A topological index is a numeric quantity of a molecule that is mathematically derived in unambiguous
way from the structural graph of a molecule. In theoretical chemistry, distance-based molecular structure descriptors are used for
modeling physical, pharmacologic, biological and other properties of chemical compounds.
Arguably the best known of these indices is the Wiener index $W(G)$, defined as the sum of all distances between distinct vertices,
$$W(G) = \frac{1}{2} \sum_{u \in V} \sum_{v \in V} d_G (u, v).$$

The $n$-cube $Q_n$ is the Cartesian product of $n$ factors $K_2$, that is $Q_n = \square_{i = 1}^n K_2$. It can be equivalently described as the graph whose vertex set consists of all $n$-tuples $b_1b_2\ldots b_n$ with $b_i \in \{0, 1\}$, where two vertices are adjacent if the corresponding tuples differ in precisely one position \cite{ImKl00}.

A subgraph $H$ of a graph $G$ is isometric if for any vertices $u$ and $v$ of $H$, $d_H(u, v) = d_G(u, v)$.
The class of graphs that consists of all isometric subgraphs of hypercubes is called partial cubes. We point out that hypercubes, even cycles, trees, median graphs (in particular acyclic cubical complexes), benzenoid graphs,
phenylenes, and Cartesian products of partial cubes are all partial cubes. A subgraph of a graph is called convex if for any two vertices of the subgraph all shortest paths (of the entire graph) between them belong to the subgraph.

Let $G$ be a connected graph. Then $e = xy$ and $f = uv$ are in the Djokovi\' c-Winkler relation
iff
$$d_G(x, u) + d_G(y, v) \neq d_G(x, v) + d_G(y, u).$$
The relation $\Theta$ is always reflexive and symmetric, and is transitive on partial cubes. Therefore,
$\Theta$ partitions the edge set of a partial cube $G$ into equivalence classes $F_1, F_2, \ldots, F_k$, called $\Theta$-classes (or cuts).

A mathematical result is really valuable only if it has been discovered independently at least twice. Note that the definition of the relation $\Theta$ is due to Djokovi\' c \cite{Dj73}, and Winkler \cite{Wi84} defined another relation which in bipartite graphs coincides with this definition -- hence it is commonly called Djokovi\' c-Winkler relation.

\section{Cut method}

The cut method was initiated in \cite{KlGuMo95} by Klav\v zar, Gutman and Mohar, where it was shown how
cuts can be used to compute the Wiener index of graphs which admit isometric embeddings into hypercubes.
One of the first applications was proving the closed formula for the Wiener index of circumcoronene series $H_k$ and solving elegantly this open problem.
The result was extended in \cite{Kl06} to general graphs by establishing a connection between the Wiener index of a graph and its canonical metric representation. For a recent survey of the cut method and overview of its applications in chemical graph theory see \cite{Kl08}.

\begin{theorem}[Cut method] \label{sandi}
Let $G$ be a partial cube and let $F_1, \ldots, F_k$ be its $\Theta$ classes. Let $n_1(F_i)$ and
$n_2(F_i)$ be the number of vertices in the two connected components of $G \setminus F_i$. Then
$$W (G) = \sum_{i = 1}^k n_1 (F_i) \cdot n_2 (F_i).$$
\end{theorem}

Our primary motivation was the recent paper \cite{MaRa13} by Manuel, Rajasingh, Rajan and Rajan, in which the cut method was reinvented and applied to compute the Wiener index of certain hexagonal and octagonal nano structures such as $C_4C_8(S)$ nanosheet, $H$-naphtalenic nanosheet and pericondensed benzenoid graphs. The authors apparently missed the whole series of papers that deal with computing distance based indices using the cut method. The same method was successively applied to hyper Wiener index \cite{CaKlPe02}, Szeged index \cite{GuKl95}, PI index \cite{Kl07}, the edge-Wiener and edge-Szeged index \cite{YiKhAs11}, degree distance \cite{IlKlSt10}, weighted Wiener index \cite{KlNa13}, terminal Wiener index \cite{IlIl13}, etc. For some other applications and properties see \cite{ChKl97}, \cite{DiCiJo08},  \cite{KlGu97}, \cite{KlGu03}.

Also the cut method was recently rediscovered by Yousefi-Azari, Khalifeh and Ashrafi in \cite{YiKhAs11} for the edge versions of the Wiener and Szeged indices, as pointed out in \cite{KlNa13}.

The main result from \cite{MaRa13} is the following.

\begin{theorem}[I-Partition Lemma] Let $G$ be a graph on $n$ vertices. Let $\{S_1, S_2, \ldots, S_m\}$ be a
partition of $E(G)$ such that each $S_i$ is an edge cut of $G$ and the removal of edges of $S_i$ leaves $G$ into
two components $G_i$ and $G_i'$. Also each $S_i$ satisfies the following conditions:
\begin{enumerate}[($i$)]
\item For any two vertices $u, v \in G_i$, a shortest path between $u$ and $v$ has no edges in $S_i$.
\item For any two vertices $u, v \in G_i'$, a shortest path between $u$ and $v$ has no edges in $S_i$.
\item For any two vertices $u \in G_i$ and $v \in G_i'$, a shortest path between $u$ and $v$ has exactly one edge in $S_i$.
\end{enumerate}

Then $W(G) = \sum_{i = 1}^m |V (G_i)| \cdot |V (G_i')|$.\label{ipartition}
\end{theorem}

Note that conditions $(i)$ and $(ii)$ are equivalent to $G_i$ and $G_i'$ being convex subgraphs. The condition $(iii)$ is in fact redundant, and it follows from the first two conditions. Namely, the shortest path from $u \in G_i$ and $v \in G_i'$ must have an odd number of edges from the cut $S_i$. Assume that the shortest path between the vertices $u \in G_i$ and $v \in G_i'$, has at least three edges in $S_i$, then one can find a shorter path from $u$ in $G_i$ to some of the other vertices from $G_i$ -- which contradicts the fact that $G_i$ and $G_i'$ are convex.

Using the following proposition from \cite{KlNa13}, the equivalence between Theorem \ref{ipartition} and Theorem \ref{sandi} is obvious:

\begin{proposition}
Let $G$ be a connected graph. Then $G$ admits a partition $\{F_i\}$ of
$E(G)$ such that $G \setminus F_i$ is a two component graphs with convex components if and
only if $G$ is a partial cube.
\end{proposition}

Also note that the formula for Wiener index for the $H$-naphtalenic nanosheet ($2n$, $2m$) from Theorem 3.3 in \cite{MaRa13} is not completely accurate (the case $2n < m$ is missing) and it is not given in the closed form (the summation cannot depend on the iterator $j$). However the formula seems to have very complicated closed form and hence not worth further research.

\section{Further generalization}

Chepoi, Deza and Grishukhin \cite{ChDeGr97} extended Theorem \ref{sandi} from the partial cubes to the class of all $L_1$-graphs that contains also many chemical non-bipartite graphs. A graph $G$ is an $L_1$-graph if it admits a scale embedding into a hypercube, where a scale embedding of $H$ into $G$ is a mapping $\beta : V(H) \rightarrow V (G)$ such that $d_G(\beta(u), \beta(v)) = \lambda \cdot d_H(u, v)$ holds for some fixed integer $\lambda$ and all vertices $u, v \in V(H)$. Hence a scale embedding with $\lambda = 1$ is an isometric embedding. 

The following result from \cite{MaRa13} is exactly Proposition 5 from \cite{ChDeGr97}.

\begin{theorem}[kI-Partition Lemma] Let $G$ be a graph on $n$ vertices. Let $E^k(G)$ denote a
collection of edges of $G$ with each edge in $G$ repeated exactly $k$ times.  Let $\{S_1, S_2, \ldots, S_m\}$ be a
partition of $E^k(G)$ such that each $S_i$ is an edge cut of $G$ and the removal of edges of $S_i$ leaves $G$ into
two components $G_i$ and $G_i'$. Also each $S_i$ satisfies the following conditions:
\begin{enumerate}[($i$)]
\item For any two vertices $u, v \in G_i$, a shortest path between $u$ and $v$ has no edges in $S_i$.
\item For any two vertices $u, v \in G_i'$, a shortest path between $u$ and $v$ has no edges in $S_i$.
\item For any two vertices $u \in G_i$ and $v \in G_i'$, a shortest path between $u$ and $v$ has exactly one edge in $S_i$.
\end{enumerate}

Then $W(G) = \frac{1}{k} \sum_{i = 1}^m |V (G_i)| \cdot |V (G_i')|$.
\end{theorem}

\begin{proposition}
Let $G$ be a finite graph scale $\lambda$ embeddable into a hypercube and let $F_1, \ldots, F_k$ be the family of (not necessarily distinct) convex cuts defining this embedding. Then, 
$W (G) = \frac{1}{\lambda}\sum_{i = 1}^k n_1 (F_i) \cdot n_2 (F_i)$.
\end{proposition}

\vspace{0.5cm} {\bf Acknowledgment. } The author is grateful to Sandi Klav\v zar for 
his remarks and discussions that helped to improve the article.

\end{document}